\documentclass{article}
\usepackage{graphicx} 

\usepackage{lineno}

\usepackage{color}

\usepackage[cmex10]{amsmath}
\usepackage{amsfonts, bm}

\usepackage{enumerate}
\usepackage{url}
\usepackage{color,hyperref}
\usepackage{tabularx}
\usepackage{amsthm}
\usepackage{amssymb,amsfonts,textcomp}
\usepackage{latexsym}
\usepackage{graphicx} 
\usepackage{multirow}
\usepackage{microtype}
\usepackage{rotating}
\usepackage{url}
\usepackage{color}
\usepackage{comment}
\usepackage{paralist}
\usepackage{booktabs}
\usepackage{bm}
\usepackage[]{units}
\usepackage{tikz}
\usepackage{lmodern}
\usepackage{tcolorbox,complexity,multirow}
\usepackage{float}

\usepackage[ruled,vlined]{algorithm2e}
\usepackage{algorithmic}

\usepackage{tikz}

\usetikzlibrary{arrows,decorations.markings,patterns,calc,positioning,backgrounds,bending}

\tikzstyle{vertex} = [circle, draw=black, scale= 0.7]
\tikzstyle{edgelabel} = [circle, fill=white, scale= 0.9]
\tikzstyle{dashedpointedline} = [line width=0.2mm,dashed,dash pattern=on 2mm off 1mm,
decoration={markings,
	mark=at position 1 with {\arrow[line width=0.2mm, scale= 1.8]{>}}
},
postaction={decorate}
] 
\tikzstyle{pointedline} = [line width=0.3mm,
decoration={markings,
	mark=at position 1 with {\arrow[line width=0.2mm, scale= 2]{>}}
},
postaction={decorate}
] 

\newenvironment{problem}[1][htb]
  {
   \begin{algorithm}[#1]%
  }{\end{algorithm}}
  \SetKwInOut{Input}{Input}
\SetKwInOut{Output}{Output}

\usepackage{todonotes}
 \def\showtodonotes{1}

\ifx\showtodonotes\undefined
    \newcommand{\erika}[1]{}
\else
    \newcommand{\erika}[1]{\todo[linecolor=orange,backgroundcolor=orange!35,bordercolor=orange,size=\tiny]{Erika: #1}}

\newcommand{\add}{\color{black}}

\newcommand{\finbox}{\hspace*{\fill}$\rule{0.17cm}{0.17cm}$}
\newcommand{\FBOX}{\hspace*{\fill}$\rule{0.17cm}{0.17cm}$}
\newcommand{\BB}{\hspace*{\fill}$\rule{0.17cm}{0.17cm}\
\rule{0.17cm}{0.17cm}$}

\newcommand{\finboxHere}{\ $\rule{0.17cm}{0.17cm}$}
\newcommand{\BOX}{\ $\rule{0.17cm}{0.17cm}$}

\newcommand{\ol}{\overline}

\newcommand{\ola}{\stackrel{\leftarrow}}
\newcommand{\ora}{\stackrel{\rightarrow}}

\newcommand{\olra}{\stackrel{\leftrightarrow}}

\def\eref#1{(\ref{#1})}

\def\emeref#1{ {\em(\ref{#1})} } \def\emeref#1{ {\em(\ref{#1})}}

\newcommand{\Proof}{\noindent {\bf Proof.  }}

\newtheorem{theorem}{Theorem}[section] \newtheorem{lemma}[theorem]{Lemma}
\newtheorem{claim}[theorem]{Claim} \newtheorem{COR}[theorem]{Corollary}
\newtheorem{PROP}[theorem]{Proposition}

\newtheorem{Thm}[theorem]{Theorem} \newcommand{\Theorem}{\begin{Thm}}
\newcommand{\eTh}{\end{Thm}}

\newtheorem{CON}[theorem]{Conjecture}
\newcommand{\Conjecture}{\begin{CON}} \newcommand{\eCon}{\end{CON}}

\newtheorem{CONS}[theorem]{} \newcommand{\Conjectures}{\begin{CONS}}
\newcommand{\eCons}{\end{CONS}}

\newcommand{\eq}{\begin{equation}} \newcommand{\eeq}{\end{equation}}

\newcommand{\THEOREM}{\begin{theorem}} \newcommand{\eT}{\end{theorem}}
\newcommand{\Lemma}{\begin{LEMMA}} \newcommand{\eL}{\end{LEMMA}}
\newcommand{\Claim}{\begin{CLAIM}} \newcommand{\eCl}{\end{CLAIM}}
\newcommand{\Corollary}{\begin{COR}} \newcommand{\eCo}{\end{COR}}

\newcommand{\Proposition}{\begin{PROP}} \newcommand{\eP}{\end{PROP}}

\newtheorem{EXAM}{Example}[section]

\newcommand{\Example}{\begin{EXAM}} \newcommand{\eXa}{\end{EXAM}}

\newtheorem{EX}{Exercise}[section]

\newcommand{\Exercise}{\begin{EX}} \newcommand{\eEx}{\end{EX}}
\newtheorem{EXS}[EX]{} \newcommand{\Exercises}{\begin{EXS}}
\newcommand{\eExs}{\end{EXS}}

\newtheorem{PROBLEM}[EX]{Problem}
\newcommand{\Problem}{\begin{PROBLEM}}
\newcommand{\ePm}{\end{PROBLEM}}

\newtheorem{PROB}[EX]{} \newcommand{\Prob}{\begin{PROB}}
\newcommand{\eProb}{\end{PROB}}

\newtheorem{PROBLEMS}[EX]{} \newcommand{\Problems}{\begin{PROBLEMS}}
\newcommand{\ePms}{\end{PROBLEMS}}

\newtheorem{OPROBLEM}[EX]{Research problem}
\newcommand{\Open}{\begin{OPROBLEM}} 
\newcommand{\eO}{\end{OPROBLEM}}

\newtheorem{REMARK}[theorem]{Remark}
\newcommand{\Remark}{\begin{REMARK}} 
\newcommand{\eRe}{\end{REMARK}}

\newcommand{\eCom}{\end{COMMENT}}

\newtheorem{TODO}{Todo}[section] \newcommand{\Todo}{\begin{TODO}}
\newcommand{\eTo}{\end{TODO}}

\newtheorem{MEMO}{Memo}[section] \newcommand{\Memo}{\begin{MEMO}}
\newcommand{\eMe}{\end{MEMO}}

\title{A network flow approach to \break a common generalization of
\break Clar and Fries numbers}
\author{Erika Bérczi-Kovács  \thanks{Alfr\'ed R\'enyi Institute of Mathematics, Department of Operations Research, ELTE E\"{o}tv\"{o}s Lor\'and University, Budapest, Hungary, and ELKH-ELTE Egerv\'ary Research Group on Combinatorial Optimization. Contact: erika.berczi-kovacs@ttk.elte.hu }
\and András Frank \thanks{Department of Operations Research, ELTE E\"{o}tv\"{o}s Lor\'and University, Budapest, Hungary, and ELKH-ELTE Egerv\'ary Research Group on Combinatorial Optimization. Contact: andras.frank@ttk.elte.hu}}
\date{\today}

\begin{document}

\maketitle

\begin{abstract}
Clar  number  and  Fries  number  are  two  thoroughly  investigated  parameters of plane graphs emerging from mathematical chemistry to measure stability of organic molecules.  
First, we introduce a common generalization of these two concepts for bipartite plane graphs, and then we extend it further to the notion of source-sink pairs of subsets of
nodes in a general (not necessarily planar) directed graph. 
The main result is a min-max formula for the maximum weight of a source-sink pair.  The proof is based on the recognition that the convex hull of source-sink pairs can be obtained as  the  projection  of  a  network  tension  polyhedron.   The  construction  makes  it  possible  to apply any standard cheapest network flow algorithm to compute both a maximum weight source-sink pair and a minimizer of the dual optimization problem formulated in the min-max theorem. As a consequence,  our  approach  gives  rise  to the first purely combinatorial, strongly polynomial algorithm to  compute  a  largest  (or  even  a  maximum  weight)  Fries-set of a perfectly matchable plane bipartite graph and an optimal solution to the dual minimization problem.\end{abstract}

\section{Introduction}

This paper considers some optimization problems arising in
mathematical chemistry.  A natural graph model to capture the basic structure of an organic
molecule uses a planar graph $G$ where the atoms are represented by
nodes of $G$ while the bonds bet\-ween atoms correspond to the edges
of $G$.  
Several important parameters can be described with this approach,
which gives rise to challenging problems in combinatorial
optimization.  One of the most
well-known of these problems concerns the Clar number
\cite{Clar-book}.  Clar's empirical observation  states a connection between the
stability of a so-called aromatic hydrocarbon molecule and the
maximum number of disjoint benzenoid rings, called the Clar number.
In the graph model, the Clar number of a molecule is the cardinality
of the largest Clar-set, where a Clar-set of a perfectly matchable
planar graph $G$ embedded in the plane consists of some disjoint inner faces of $G$ that can simultaneously alternate
with respect to an appropriate perfect matching of $G$.

Another important parameter of a molecule is its Fries number
\cite{Fries, Hansen-Zheng94b}.  In the graph model, the Fries number of a perfectly
matchable planar graph $G$ embedded in the plane is the maximum
cardinality of a Fries-set, where a Fries-set consists of some
distinct but not necessarily disjoint inner faces that can
simultaneously alternate with respect to an appropriate perfect
matching of $G$.


A major issue of the investigations is the algorithmic and structural
aspects of Clar and Fries numbers and sets, in which the polyhedral
approach (that is, linear programming) plays a central role.  We
concentrate on 2-connected perfectly matchable bipartite plane graphs
$G=(S,T;E)$, where a plane graph means a planar graph with a given
embedding in the plane.  (Recall that a polyhedron is the solution set of a linear inequality system, while a polytope is the convex hull of a finite number of points.  A fundamental theorem of linear programming states that every polytope is a bounded polyhedron.  In discrete optimization, it is often a basic, non-trivial approach to provide explicitly such a linear description of the polytope of combinatorial objects to be investigated.)

Consider first the situation with Clar numbers and sets.  Hansen and
Zheng \cite{Hansen-Zheng94a} (using different notation) introduced a
polyhedron $P=\{(x',x''):  Fx'+Qx''=1, \ (x',x'')\geq 0\}$, where $F$ is
the node vs inner-face incidence matrix of $G$ while $Q$ is the
standard node-edge incidence matrix of $G$ (and hence $[F,Q]$ is a
small size $(0,1)$-matrix).  They considered the projection $P_{\rm
Cl} := \{x':  \ \exists \ x'' $ \ with \ $(x',x'')\in P\}$ of $P$ to
the variables $x'$ and observed that the integral (that is,
$(0,1)$-valued) elements of $P_{\rm Cl}$ correspond to the Clar-sets
of $G$.  Hansen and Zheng conjectured that $P$ (and hence $P_{\rm
Cl}$) is an integral polyhedron, which property, if true, implies that
the vertices of $P_{\rm Cl}$ correspond to the Clar-sets of $G$.

Abeledo and Atkinson \cite{Abeledo-Atkinson3} (Theorem 3.5) proved
that $[F,Q]$ is a unimodular matrix (though not always totally
unimodular).  This implied the truth of the conjecture of Hansen and
Zheng since a theorem of Hoffman and Kruskal \cite{Hoffman-Kruskal}
(see Theorem 21.5 in the book \cite{Schrijverbook0} of Schrijver)
implies that, for a unimodular matrix $[A,A']$ and integral vector
$b$, the polyhedron $\{(x,x'):  Ax+A'x'=b, x\geq 0\}$ is integral (and
so is the dual polyhedron $\{y:  yA\geq c, yA'=c' \}$ whenever
$(c,c')$ is an integral vector.)

A tiny remark is in order.  It follows from the previous approach that
$P_{\rm Cl}$ is the polytope of Clar-sets, and $P_{\rm Cl}$ was
obtained as a projection of the polyhedron $P$.  At first sight, it
may seem more natural to figure out an explicit polyhedral description
of $P_{\rm Cl}$ itself.  However, the linear inequality system
describing $P_{\rm Cl}$ may consist of an exponential number of
inequalities.  This is why it is more efficient to produce $P_{\rm
Cl}$ as the projection of the (integral) polyhedron $P$ that is
described by only $O(\vert E\vert )$ constraints.  Of course, a
maximum weight vertex of $P_{\rm Cl}$ (that is, a maximum weight
Clar-set) can be computed once a maximum weight vertex of $P$ can be.

Abeledo and Atkinson \cite{Abeledo-Atkinson3} noted that an analogous
approach works for Fries-sets, as well.  Namely, in Theorem 3.6 of
\cite{Abeledo-Atkinson3} they claimed (without a proof) that a certain
matrix $A$ is unimodular, and indicated that the Fries-sets correspond
to the vertices of the a projection $P_{\rm Fr}$ of the polyhedron
defined by $A$.  It should be noted that $A$ is a small size $(0,\pm
1)$-matrix.

Since the sizes of the describing unimodular matrices are small (a
linear function of the size of the graph $G$), a general purpose
linear programming algorithm can be used to compute an optimal
(maximum weight) basic solution which corresponds to an optimal
Clar-set (Fries-set).  Moreover, as the matrices are $(0,\pm
1)$-valued, the strongly polynomial time algorithm of Tardos
\cite{Tardos86} for solving such linear programs can be applied.
Therefore, there exists a strongly polynomial algorithm for computing
both a maximum weight Clar-set and a maximum weight Fries-set.

In addition, Abeledo and Atkinson \cite{Abeledo-Atkinson4} provided a
network flow formulation for the dual of the Clar problem, and hence,
as they write:  {\it "the Clar number and the cut cover can be found
in polynomial time by solving a minimum cost network flow problem"}.
Yet another fundamental result of Abeledo and Atkinson
\cite{Abeledo-Atkinson4} is an elegant min-max formula for the Clar
number of a plane bipartite graph, which was conjectured by Hansen and
Zheng \cite{Hansen-Zheng92,Hansen-Zheng94a}.

As for a min-max formula concerning the Fries number, the
above-mentioned linear programming approach of Fries-sets
\cite{Abeledo-Atkinson3} does provide one through the l.p. duality
theorem, though this was not explicitly formulated (in a form similar
to the one mentioned before concerning Clar numbers).  Moreover, to
our best knowledge, no purely combinatorial or network flow based
algorithm appears in the literature for computing the Fries number of
$G$.

The major goal of the present work is to develop a common
generalization of the Clar number and Fries number problems, and based
on this, to develop a strongly polynomial algorithm relying
exclusively on standard network flow subroutines.  Our approach is
self-contained, in particular, it does not rely on any earlier
polyhedral results on Clar-sets and Fries-sets.

We note that the presented common framework for Clar-sets and
Fries-sets has nothing to do with a popular conjecture stating that a
plane graph has a largest Clar-set that can be extended to a largest
Fries-set.  This conjecture was disproved by Hartung \cite{Hartung14}
for general plane graphs (for fullerenes), proved for a special class of benzenoids by
Graver et al.  \cite{GHS13}, but it is still open for general
perfectly matchable bipartite plane graphs.


A basic tool of our approach is that we reformulate the original
problem on $G$ for the planar dual graph of $G$, and then we can drop
the planarity assumption.  This idea appeared in the work of Erd{\H
o}s, Frank, and Kun \cite{FrankJ66}, where the major problem was to
find a largest (weight) sink-stable set of a (not necessarily planar)
digraph $D=(V,A)$. Here the term sink-stable refers to a stable set of $D$ whose elements are sink nodes in an appropriate dicut-equivalent reorientation of $D$, where a digraph $D'=(V,A')$ is a dicut-equivalent reorientation of $D$ if it can be obtained from $D$ by reorienting some disjoint one-way cuts ($=$ dicuts).
(Incidentally we remark that this digraph model not only covered the
(weighted) Clar number problem but implied an elegant theorem of Bessy
and Thomass\'e \cite{Bessy-Thomasse} as well as its extension by
Seb{\H o} \cite{Sebo07}.)

The Fries number problem can also be embedded into this general
digraph framework, where it transforms into finding a dicut-equivalent
reorientation of a digraph in which the set of sink and source nodes
is as large as possible.
However, this requires a technique more involved than the one needed for Clar numbers (see Remark \ref{kornemeleg}).

Our goal, in fact, is to develop a common framework for these extensions of the Clar and Fries problems.  Namely, in Section \ref{subsec:notation}, we introduce the notion of source-sink (so-si) pairs
$\{Y_{\text{o}}, Y_{\text{i}}\}$ of disjoint subsets of nodes of $D$, and, for a given pair of weight-functions $\{w_{\text{o}}, w_{\text{i}}\}$, develop a min-max formula (Theorem \ref{minmax2})  for the maximum total weight $\widetilde{w}_{\text{o}}(Y_{\text{o}}) +\widetilde{w}_{\text{i}}(Y_{\text{i}})$ of a so-si pair. The  proof  is  based  on  the  revelation  that  the  convex  hull  of  (the  characteristic  vectors of) so-si pairs can be obtained as the projection of a network tension polyhedron (of $3|V|$-dimension),  and hence the maximum weight so-si pair problem can be transformed into a maximum weight feasible tension problem.  As such a tension problem is the linear programming dual of a minimum cost network flow (or circulation) problem, the approach not only results in a (combinatorial) min-max formula but it also gives rise to a strongly polynomial algorithm that relies exclusively on standard network flow subroutines.

\subsection{Notation, terminology}\label{subsec:notation}

Given a ground-set $S$, we do not distinguish between a vector $c\in
{\bf R} \sp S$ and a function $c:S\rightarrow {\bf R}$.  For a vector
or function $c$ on $S$, let $\widetilde c(X):=\sum \ [c(s):s\in X]$.

We call a digraph {\bf weakly connected} or just {\bf connected} if
its underlying undirected graph is connected.  Often we refer to the
edges of a digraph as {\bf arc}s.  Let $D=(V,A)$ be a loopless
connected digraph with $n\geq 2$ nodes.  For a subset $U$ of nodes,
the {\bf in-degree} $\varrho (U)=\varrho _D(U)=\varrho _A(U)$ is the
number of arcs entering $U$, while the {\bf out-degree} $\delta
(U)=\delta _D(U)=\delta _A(U)$ is the number of arcs leaving $U$.
Typically, we do not distinguish between a one-element set $U=\{u\}$
and its only element $u$, for example, $\varrho (u)=\varrho (\{u\})$.
For a vector (or function) $z$ on the arc-set $A$ of $D$, let $\varrho
_z(U):= \sum \ [z(e):  e$ enters $U]$ and $\delta _z(U):= \sum \
[z(e):  e$ leaves $U]$.

A function $z$ defined on the arc-set of $D$ is called a {\bf
circulation} in $D$ if $\varrho _z(v) = \delta _z(v)$ holds for every node
$v$.  This is equivalent to requiring that $\varrho _z(v) \leq \delta
_z(v)$ holds for every node $v$.  Furthermore, for a circulation $z$,
$\varrho _z(U)=\delta _z(U)$ holds for any subset $U\subseteq V$.  For
lower and upper bound functions $f,g$, we say that a circulation $z$
is {\bf feasible} if $f\leq z\leq g$.  (Here $-\infty $ is allowed for
the components of $f$ and $+\infty $ for the components of $g$.)

A function $\pi $ defined on the node-set of $D$ is often referred to
a {\bf potential}, which {\bf induces} the {\bf tension} $\Delta _\pi
$ (sometimes called potential-drop) on the arc-set of $D$, where
$\Delta _\pi (uv):= \pi (v)-\pi (u)$ for $uv\in A$.  An integer-valued
tension can be obtained as the potential-drop of an integer-valued
potential.  It is a well-known property that the set of circulations
and the set of tensions are complementary orthogonal subspaces (over
the rationals or the reals), that is, $z\Delta _\pi =0$ holds for any
circulation $z$ and potential $\pi $. We call an integer-valued
potential $\pi $ {\bf small-dropping} if its potential-drop $\Delta
_\pi $ is $(0,1)$-valued.  When $\Delta _\pi $ is $(0,1)$-valued only
on a subset $A_0\subseteq A$ of arcs, we say that $\pi $ is {\bf
small-dropping on} $A_0$.

For given lower bound $c_\ell:A\rightarrow {\bf R} \cup \{-\infty \}$
and upper bound $c_u:A\rightarrow {\bf R} \cup \{+\infty \}$, we say
that a potential $\pi $ or the potential drop $\Delta _\pi $ is {\bf
$(c_\ell, c_u)$-feasible} if $c_\ell \leq \Delta _\pi \leq c_u$.  When
$c_\ell\equiv -\infty $ (that is, when no lower bound is given), we
speak of {\bf $c_u$-feasibility}.

A {\bf circuit} is a connected undirected graph in which the degree of
every node is 2. For a circuit $C$, we use the convention that $C$
also denotes the edge-set of the circuit while $V(C)$ denotes its
node-set.  A directed graph is also called a {\bf circuit} if it
arises from an undirected circuit by arbitrarily orienting its edges.
A directed graph obtained from an undirected circuit by orienting each
edge in the same direction is called a {\bf one-way circuit} or just a
{\bf di-circuit}.

We call two orientations of an undirected graph {\bf di-circuit
equivalent} if one of them can be obtained from the other one by
reorienting a set of arc-disjoint di-circuits.  This is clearly
equivalent to requiring that the in-degree of a node on one
orientation is equal to the in-degree in the other orientation.

By {\bf reorienting} (or {\bf reversing}) an arc $uv$, we mean the
operation of replacing $uv$ by $vu$.  The reorientation of a subset
$B$ of arcs (that is, reversing $B$) means that we reorient all the
elements of $B$.

By a {\bf cut} of a connected graph $G=(V,E)$, we mean the set of
edges connecting $Z$ and $V-Z$ for some subset $Z$ of nodes. 
By a \textbf{minimal} cut we mean an inclusionwise minimal cut.
A useful observation is that a cut is minimal if and only if both $Z$ and $V-Z$ induce a connected subgraph of $G$.  Furthermore, every cut can
be partitioned into minimal cuts.  An minimal cut is sometimes
called a {\bf bond}, but this term should not be confused with a bond
used in chemistry (which corresponds to an edge of the modelling
graph).

By a {\bf cut} of a digraph $D=(V,A)$ defined by a subset $Z$ of
nodes, we mean the set of arcs connecting $Z$ and $V-Z$ (in either
direction).  In the special case when no arc enters $Z$, the cut is
called a {\bf directed cut} or a {\bf dicut} of $D$.  Sometimes the
term {\bf one-way cut} is used for a dicut.  
We say that digraph $D'= (V, A')$ is called a \textbf{dicut-equivalent reorientation} of $D$ if $D'$ is obtained from $D$ by reorienting some disjoint one-way cuts of $D$ (see \cite{FrankJ66}).
A node of $D$ will be
called a {\bf sink node} (or just a sink) if it admits no leaving
arcs.  A node is a {\bf source node} if it admits no entering arcs.  A
subset of nodes is a {\bf sink set} (respectively, a {\bf source set})
if each of its elements is a sink (resp., a source).  Clearly, a sink
set is always stable.

We call a subset $Y$ of nodes of $D$ {\bf sink-stable} ({\bf
source-stable}) if there is a dicut-equivalent reorientation of $D$ in
which every node in $Y$ is a sink node (source node).  Obviously, a
set is sink-stable if and only if it is source-stable.  These are
special stable sets of $D$, and a node of a one-way circuit of $D$ never
belongs to a sink-stable set of $D$.

For two disjoint sets $Y_{\rm o}, Y_{\rm i}\subseteq V$, we call the
ordered pair $(Y_{\rm o},Y_{\rm i})$ a {\bf source-sink} (= {\bf
so-si}) {\bf pair} if there is a dicut-equivalent reorientation of $D$
in which $Y_{\rm o}$ is a source set while $Y_{\rm i}$ is a sink set.
It may be the case that $(Y_{\rm o},Y_{\rm i})$ is a so-si pair, but
$(Y_{\rm i},Y_{\rm o})$ is not (as exemplified by an acyclic
orientation of a triangle).  For a subset $Y$ of nodes the following four properties
are obviously equivalent:  $Y$ is a sink-stable set, $Y$ is a source-stable set,
$(Y,\emptyset )$ is a so-si pair, $(\emptyset ,Y)$ is a so-si pair.

We call a subset $Y$ of nodes of $D$ {\bf resonant} if it can be
partitioned into subsets $Y_{\rm o}$ and $Y_{\rm i}$ where $(Y_{\rm
o},Y_{\rm i})$ is a so-si pair.  (The name is motivated by the term
resonant set of faces of a plane graph, see, for example,
\cite{ZAY18}).  A sink-stable (or source-stable) set is clearly
resonant\footnote{The term 'resonant' has different definitions in the literature, sometimes it is a synonym of a Clar-set. In this paper we use it as a generalization of a Fries-set.}.

\medskip



\section{Undirected plane graphs and directed general graphs}

Recall a basic result of F\'ary \cite{Fary} stating that every simple
planar graph admits an embedding in the plane where each edge is
represented by a straight-line segment.  A deeper result of Tutte
\cite{Tutte63} states that a 3-connected simple planar graph admits an
embedding into the plane in which every edge is represented by a
straight-line segment and the faces are convex sets.

Let $G=(V,E)$ be a loopless 2-connected planar graph.  Throughout we
consider a specified embedding of $G$ into the plane, and say that $G$
is a {\bf plane} graph.  Such an embedded graph divides the plane into
disjoint (connected) regions, exactly one of them is infinite.  We
refer to these regions as the {\bf face}s of the plane graph.  Since
$G$ is 2-connected, its faces are surrounded by a circuit of $G$, and
we refer to these as {\bf face-bounding circuits} or just {\bf
face-circuits}.  (In the literature, often the longer term boundary
circuit is used.)  With a slight abuse of terminology, sometimes we do
not distinguish between a face and its face-circuit.  The single
infinite face is called the {\bf outer face} of $G$ while all the
other faces are the {\bf inner} faces of $G$.  For their surrounding
circuits, we use the term {\bf outer} ({\bf inner}) face-circuit.

Note that the set of face-circuits may be different for another
embedding of a planar graph into the plane, but for 3-connected planar
graphs, by a theorem of Tutte, the set of face-circuits is independent from the embedding.

Let $F$ be a face of $G$ and $C_F$ its face-circuit.  Let $\ora G$ be
an orientation of $G$ in which the orientation $\ora C_F$ of $C_F$ is
a one-way circuit.  When we say that the face $F$ and the di-circuit
$\ora C_F$ is clockwise (anti-clockwise) oriented, we mean that it is
clockwise (anti-clockwise) when looked at from the inside of $F$.  (Or
putting it in another way, we may imagine that $\ora G$ is embedded on
the surface of a sphere, with no distinction between the outer face
and the inner faces, and a one-way face-circuit is clockwise in the
usual sense.)

As a simple example, let the graph be a single circuit $C$ embedded in
the plane and let $\ora C$ be an orientation of $C$ which is a one-way
circuit.  Now the graph has one inner face and one outer face, $C$ is
a one-way face-circuit belonging to both faces, and $\ora C$ is
clockwise with respect to the inner face if and only if it is
anti-clockwise with respect to the outer face.

It follows from this definition that if $F_1$ and $F_2$ are
neighbouring (incident) faces of $G$ (that is, they have an edge in
common) and their face-circuits are one-way circuits in $\ora G$, then
one of $F_1$ and $F_2$ is clockwise in $\ora G$ while the other is
anti-clockwise.

Throughout the paper, we assume that $G=(S,T;E)$ is a perfectly
matchable, 2-connected, loopless, bipartite plane graph.  Let $M$ be a
perfect matching of $G$. A circuit $C$ of $G$ is called {\bf
$M$-alternating} 
if edges in $C\cap M$ and $C- M$ are alternating on $C$. 
If $M'$ is another perfect matching of $G$, then the symmetric difference
$M\ominus M':=(M- M') \cup (M'- M)$ of $M$ and $M'$ consists of disjoint circuits which are
both $M$-alternating and $M'$-alternating.

A set ${\cal C}'$ of node-disjoint face-circuits of $G$ forms a {\bf
Clar-set} if $G$ has a perfect matching $M$ such that each member of
${\cal C}'$ is $M$-alternating.  

A set ${\cal C}'$ of (distinct)
face-circuits of $G$ forms a {\bf Fries-set} if $G$ has perfect
matching $M$ such that each member of ${\cal C}'$ is $M$-alternating.
Sometimes a Fries-set is called a resonant set of faces or
face-circuits \cite{ZAY18}.

Let ${\cal C}_G$ denote the set of face-circuits of $G$ (including the
outer face).  For a non-empty subset ${\cal C}'$ of ${\cal C}_G$, the
{\bf Clar number} Cl$({\cal C}') =$ Cl$_G({\cal C}')$ of $G$ {\bf
within} ${\cal C}'$ is the maximum cardinality of a Clar-set
consisting of some members of ${\cal C}'$.  When ${\cal C}'$ is the
set of inner faces, we simply speak of the {\bf Clar number} of $G$.
It should be noted that the Clar number of two different embeddings of
a planar graph into the plane with the same outer face may be
different (unless the graph is 3-connected).

For a non-empty subset ${\cal C}'$ of face-circuits, the {\bf Fries
number} Fr$({\cal C}') =$ Fr$_G({\cal C}')$ of $G$ {\bf within} ${\cal
C}'$ is the maximum cardinality of a Fries-set consisting of some
members of ${\cal C}'$.  When ${\cal C}$ is the set of all inner
faces, we simply speak of the {\bf Fries number} of $G$.

The basic optimization problem for Clar (respectively, Fries) numbers
is \ (A) \ to compute the Clar (Fries) number of $G$ as well as a
maximizing Clar-set (Fries-set) consisting of inner face-circuits.
In a slightly more general version, we are interested in \ (B) \
finding a maximum cardinality Clar-set (Fries-set) within a specified
subset of face-circuits.  An even more general form is \ (C) \ the
weighted Clar-set (Fries-set) problem where we are given a
non-negative weight-function $w$ on the set of faces (face-circuits),
and want to construct a Clar-set (Fries-set) of largest weight.  (When
$w$ is $(0,1)$-valued we are back at Problem (B), while Problem (A) is
obtained for the special weight-function $w$ defined to be $1$ on all
inner faces and $0$ on the outer face.)  In what follows, we consider
a weighted version of a common generalization of Clar-sets and
Fries-sets.


\subsection{Common framework for Clar and Fries}

For any perfect matching $M$ of $G$, we define an orientation $\ora
G_M$ of $G$, as follows.  Orient each element of $M$ toward $S$ while
all the other edges of $G$ toward $T$.  Then a face (and its
face-circuit) of $G$ is $M$-alternating if and only if it is a one-way
circuit in the orientation $\ora G_M$.  Note that every one-way
circuit of $\ora G_M$ is clockwise or anti-clockwise.

We will say that an $M$-alternating face-circuit of $G$ is {\bf
clockwise} ({\bf anti-clockwise}) if the corresponding one-way circuit
in $\ora G_M$ is clockwise (anti-clockwise).  It follows for two
adjacent $M$-alternating faces that one of them is clockwise, while
the other one is anti-clockwise.  Furthermore, both the set of
clockwise $M$-alternating face-circuits and the set of anti-clockwise
$M$-alternating face-circuits are Clar-sets, and their union is a
Fries-set.  And conversely, every Fries-set ${\cal F}'$ can be
obtained in this way.  Indeed, if $M$ is a perfect matching for which
the members of ${\cal F}'$ are $M$-alternating, then both the set of
clockwise $M$-alternating members and the set of anti-clockwise
$M$-alternating members of ${\cal F}'$ are Clar-sets.

Observe that for any Clar-set ${\cal C}'$ there is a perfect matching
$M'$ such that each member of ${\cal C}'$ is $M'$-alternating and
clockwise.  Indeed, by definition, there is a perfect matching $M$ for
which each member of ${\cal C}'$ is $M$-alternating.  Then the perfect
matching $M'$ obtained by taking the symmetric difference of $M$ and
the union of anti-clockwise members of ${\cal C}'$ will do.

Let $w_1$ and $w_2$ be two non-negative weight-functions on the set of
faces of $G$.  In the common generalization of the weighted Fries and
Clar problems, we look for a perfect matching $M$ of $G$ for which the
sum of the total $w_1$-weight of the clockwise $M$-alternating
face-circuits plus the total $w_2$-weight of the anti-clockwise
$M$-alternating face-circuits is as large as possible.  We refer to
this common framework as the {\bf double-weighted Clar-Fries} problem, see Problem \ref{def:DWCFP} below.

In the special case when $w_1:=w$ and $w_2:=0$, we are back at the
maximum $w$-weighted Clar-set problem.  In the special case when
$w_1:=w$ and $w_2:=w$, we are back at the maximum $w$-weighted
Fries-set problem.

\begin{figure}[h]
\begin{problem}[H]
\medskip
\KwIn{A perfectly matchable 2-connected plane bipartite graph $G=(S,T,E)$ with non-negative weightings $w_1$ and $w_2$ on the set of faces of $G$.
}
\bigskip
\KwOut{A perfect matching $M$ of $G$ for which the sum of the total $w_1$-weight of the clockwise $M$-alternating face-circuits plus the total $w_2$-weight of the anti-clockwise $M$-alternating face-circuits is maximum.
}
\caption{Double-weighted Clar-Fries problem}\label{def:DWCFP}
\end{problem}
\end{figure}

We shall show that Problem \ref{def:DWCFP} can be embedded in the following, significantly more general
framework concerning directed graphs.

\begin{figure}[h]
\begin{problem}[H]
\medskip
\KwIn{A loopless, weakly connected (non-necessarily planar) digraph $D=(V,A)$ with non-negative weightings $w_{\rm o}$ and $w_{\rm i}$ on $V.$
}
\bigskip
\KwOut{A dicut-equivalent reorientation of $D$ for which the sum of the total $w_{\rm o}$-weight of source nodes plus the total $w_{\rm i}$-weight of sink nodes is maximum.
}
\caption{Double-weighted source-sink pair problem}\label{def:WSSP}
\end{problem}
\end{figure}
\begin{claim}
The double-weighted Clar-Fries problem formulated in Problem \ref{def:DWCFP} is a special case of the double-weighted source-sink pair problem formulated in Problem \ref{def:WSSP}.\end{claim}
\Proof
Let $M$ be an arbitrary perfect matching of $G$, and consider the planar dual digraph
$D_M$ of  $\ora
G_M$ {\add(where each dual arc is oriented from the left to the right with respect to the primal arc)}. The weight-functions $w_1$ and $w_2$ on the set of faces of $G$ can be associated with weight-functions $w_{\rm o}$ and $w_{\rm i}$ on the set of nodes of $D_M$.
Observe that for another perfect matchings $M'$ of $G$, the orientations $\ora G_M$ and $\ora G_{M'}$ are di-circuit equivalent and hence $D_M$ and $D_{M'}$ are dicut-equivalent. Conversely, any dicut-equivalent reorientation of $D_M$ arises in this way (that is, there is a perfect matching $M'$ of $G$ for which the orientations $\ora G_M$ and $\ora G_{M'}$ are di-circuit equivalent). Therefore there is a one-to-one correspondence between the so-si pairs of $D_M$ and the set of clockwise and anticlockwise faces of $G$ belonging to a perfect matching $M'$ of $G$.
It follows that Problem \ref{def:WSSP}, when applied to
$D_M$ in place of $D$, outputs a so-si pair 
$(Y_{\rm o}, Y_{\rm i})$ 
for which its double-weight  $\widetilde{w}_{\text{o}}(Y_{\text{o}}) +\widetilde{w}_{\text{i}}(Y_{\text{i}})$
is maximum, 
and this pair $(Y_{\rm o}, Y_{\rm i})$ corresponds to a pair of clock-wise and anti-clock-wise $M'$-alternating faces with a maximum $(w_1, w_2)$ double weight. \FBOX

\medskip
    
Note that the problem of finding a dicut-equivalent reorientation of a
digraph $D$ for which the total weight of the set of sink nodes is as
large as possible is a special case of this problem and was solved by
Erd{\H o}s, Frank, and Kun \cite{FrankJ66}.  Theorem 4.1 of
\cite{FrankJ66} states that a stable set of $D$ is sink-stable if and
only if $\vert S\cap V(C)\vert \leq \eta (C)$ holds for every circuit $C$ of $D$, where $\eta (C)$
denotes the minimum of the number of forward arcs of $C$ and the
number of backward arcs of $C$.  By extending this result, we shall
derive an analogous characterization of source-sink pairs in Section
\ref{sosichar}.  Resonant sets will be characterized in Corollary
\ref{cor:U} in terms of integer-valued circulations (and not only
circuits).  Perhaps a bit surprisingly, no characterization of
resonant sets may exist that requires an inequality only for circuits,
see Remark \ref{kornemeleg}.

In another approach, clockwise and anti-clockwise M-alternating edge sets appear in the formulation of Abeledo and Ni \cite{Abeledo}, who gave an LP description for the Fries number of plane bipartite graphs using a TU matrix.

Shi and Zhang \cite{Shi-Zhang22} considered $(4,6)$-fullerenes, and gave a formula for a variation of the Fries number problem when only the number of alternating hexagonal
faces is to be maximized. Note that this problem is a special case of the Fries number within a subset ${\cal C}'$ of face-circuits defined above, by choosing ${\cal C}'$ to be the set of hexagonal faces.

We conclude this section by remarking that the problems above
concerning Clar and Fries numbers can be considered for non-bipartite
plane graphs, as well.  It was shown by B\'erczi-Kov\'acs and
Bern\'ath \cite{Kovacs-Bernath} that computing the Clar number is {\bf
NP}-hard for general plane graphs, but the problem is open if the
number of odd faces is constant.  This is the case, for example, for
fullerene graphs (a 3-connected 3-regular planar graph with hexagonal
faces and exactly six pentagonal faces).  The status of the Fries
number problem for general plane graphs is open.
Salami and Ahmadi in \cite{Salami-Ahmadi} gave both integer and
quadratic programming formulations for the Fries number of fullerenes.

\medskip

\section{Characterizing source-sink pairs} \label{sosichar}

In this section we describe a characterization of source-sink pairs of a digraph $D=(V,A)$. Recall the notion of dicut-equivalence of two orientations of an
undirected graph.  It was observed in \cite{FrankJ66} (Lemma 3.1) that
a subset $F\subseteq A$ of arcs is the union of disjoint one-way cuts
if and only if there is an integer-valued potential $\pi $ for which
$\chi_F= \Delta _\pi $, where $\chi_F\in \{0,1\}^A$ is the incidence vector of $F$. This immediately implies the following observation.

\begin{claim} \label{vekvi1}  
{\add Digraph $D$ admits a
dicut-equivalent reorientation $D'$
if }and only if there is a potential $\pi :  V \rightarrow {\bf Z}_+$ which is small-dropping on $A$ and $D'$ arises from $D$ by reversing those arcs $e$ of $D$ for which $\Delta _\pi (e) = 1$. \FBOX \end{claim}
\medskip

{\add  Let $f:A\rightarrow {\bf R} \cup \{-\infty \}$ and $g:A\rightarrow {\bf R}
\cup \{+\infty \}$ be lower and upper bound
functions on the arc-set of $D=(V,A)$, respectively,} for which $f\leq g$.  By
the extended (but equivalent) form of Gallai-lemma \cite{Gallai58},
there is an $(f,g)$-feasible tension if and only if, for every circuit
$C$ of $D$, the total $f$-value on the arcs of $C$ in one direction is
at most the total $g$-value on the arcs in the other direction.  (In
the original Gallai-lemma $f\equiv -\infty $ in which case the
condition reduces to requiring that $g$ is conservative, that is, $D$
admits no one-way circuit with negative $g$-weight.  This immediately
implies the extended form by adding the reverse $\ola e$ of each arc
$\in A$ to $D$ and define $g(\ola e):  = -f(e)$.)

Let $F$ and $R$ be two disjoint subsets of arcs of digraph $D=(V,A)$.
We consider the members of $F$ as fixed arcs while the elements of $R$
must be reversed.  Let $N:=A-(F\cup R)$
denote the complement of $F\cup R$.

\begin{claim} \label{FT-dicut} For a partition $\{F,R,N\}$ of the arc-set of
digraph $D=(V,A)$, the following are equivalent.

\medskip

\noindent {\bf (A)} \ $D$ admits a dicut-equivalent reorientation for
which every element of $R$ is reversed while the orientation of every
element of $F$ is unchanged.

\noindent {\bf (B)} \ There is a small-dropping potential $\pi $ for
which

\eq \Delta _\pi (e)= \begin{cases} 1 & \ \ \hbox{if}\ \ \ e\in R \cr 0
& \ \ \hbox{if}\ \ \ e\in F. \end{cases} \label{(pot-felt)} \eeq

\noindent {\bf (C)} \ For every circuit of $D$, the number of $R$-arcs
in one direction is at most the number of $(R\cup N)$-arcs in the
other direction.\end{claim}

\Proof The equivalence of conditions {\bf (A)} and {\bf (B)} follows
from Claim \ref{vekvi1}.  Let $$ (f(e),g(e)):= \begin{cases} (0,0) & \
\ \hbox{if}\ \ \ e\in F \cr (1,1) & \ \ \hbox{if}\ \ \ e\in R \cr
(0,1) & \ \ \hbox{if}\ \ \ e\in N. \end{cases} $$

It follows from this definition that a small-dropping potential $\pi$ meets
\eref{(pot-felt)} if and only if $\Delta _\pi $ is $(f,g)$-feasible.
For a circuit of $D$, the sum of the $f$-values of the arcs in one
direction is the number of $R$-arcs in this direction, and the sum of
the $g$-values of the arcs of $C$ in the other direction is the number
of $(R\cup N)$-arcs in that direction.  Therefore the extended
Gallai-lemma implies the equivalence of Conditions {\bf (B)} and {\bf
(C)}.  \FBOX \medskip

In digraph $D=(V,A)$, let $Y_{\rm o}\subseteq V$ and $Y_{\rm
i}\subseteq V$ be disjoint stable sets for which we want to decide
whether $(Y_{\rm o} , Y_{\rm i})$ is a so-si pair.  Of course, if no
arc enters $Y_{\rm o}$ and no arc leaves $Y_{\rm i}$, then $(Y_{\rm o}
, Y_{\rm i})$ is a so-si pair.  We refer to an arc $e\in A$ as
\textbf{incorrect} if it enters $Y_{\rm o}$ or leaves $Y_{\rm i}$, while $e$ is
considered \textbf{correct} if it leaves $Y_{\rm o}$ or enters $Y_{\rm i}$.
We refer to those members of $A$ which are neither correct nor incorrect as \textbf{neutral} arcs.

\begin{theorem}\label{so-si.char} For disjoint node-sets $Y_{\rm o}\subseteq
V$, $Y_{\rm i}\subseteq V$ of digraph $D=(V,A)$, the following
properties are equivalent.  \medskip

\noindent {\bf (A)} \ $(Y_{\rm o},Y_{\rm i})$ is a so-si pair.

\noindent {\bf (B)} \ There is a 
small-dropping potential $\pi $ on $V$ for which 

\eq \Delta _\pi (e)= \begin{cases} 1 & \ \ \hbox{if}\ \ \ e\in A \ \
\hbox{incorrect}\ \cr 0 & \ \ \hbox{if}\ \ \ e\in A \ \
\hbox{correct.}\ \end{cases} \label{(B-felt)} \eeq

\noindent {\bf (C)} \ For every circuit of $D$, the number of
incorrect arcs in one direction is at most the total number of correct
arcs and neutral arcs in the other direction.  
\end{theorem}

\Proof By definition, $(Y_{\rm o},Y_{\rm i})$ is a so-si pair
if and only if $D$ admits a dicut-equivalent reorientation $D'$ in which
$Y_{\rm o}$ is a source set and $Y_{\rm i}$ is a sink set, that is,
the incorrect arcs in $D'$ are reversed while the correct arcs are
unchanged.  Let $R$ denote the set of incorrect arcs and $F$ the set
of correct arcs.  Then the properties {\bf (A), (B)}, and {\bf (C)}
occurring in Claim \ref{FT-dicut} are just equivalent to those in the
theorem.  \FBOX

\medskip

In the rest of the paper, we are going to prove a min-max theorem for
the maximum $w$-weight of a resonant set.  Actually, we do this in a
more general form when $w_{\rm o}:V\rightarrow {\bf R}_+$ and $w_{\rm
i}:V\rightarrow {\bf R}_+$, are two given weight-functions and we are
interested in finding a so-si pair $(Y_{\rm o},Y_{\rm i})$ whose {\bf
$(w_{\rm o},w_{\rm i})$-weight} defined by $\widetilde w_{\rm
o}(Y_{\rm o})+\widetilde w_{\rm i}(Y_{\rm i})$ is maximum.  As the
problem will be formulated as a network circulation problem in a larger
digraph $D\sp *$, a standard algorithm for the latter one will compute
the optimal solutions occurring in the min-max theorem (Theorem
\ref{minmax2}).

\medskip

\section{Min-max formula and algorithm for double-weighted source-sink pair}

In this section, we exhibit first a min-max formula for the maximum
{\bf $(w_{\rm o},w_{\rm i})$-weight} of a so-si pair of digraph
$D=(V,A)$.  
As a preparation for its proof (to be worked out in
Section \ref{mainproof}), we describe here a feasible circulation problem on a larger digraph $D\sp
*$ along with its linear programming dual, a feasible tension problem.
The algorithmic aspects of the approach will be discussed in Section \ref{algasp}.

\subsection{The main theorem}\label{subsec:MT}

Let $\ola{A}$ denote the set of arcs obtained from $A$ by reversing
all the arcs in $A$.  Let $\olra{A}:= A\cup \ola{A}$ and
$\olra{D}:=(V, \olra{A})$.  That is, $\olra{D}$ is the digraph arising
from $D$ by adding the reverse of its arcs.

For a non-negative vector $z_{\rm i}$ defined on the arc-set $\olra
A$, we say that $z_{\rm i}$ is an {\bf in-cover} of $w_{\rm i}$ if
$\varrho _{z_{\rm i}}(v)\geq w_{\rm i}(v)$ holds for every node $v\in
V$ (where $\varrho _{z_{\rm i}}(v)$ denotes the sum of $z_{i}$-values
on the arcs in $\olra{A}$ entering $v$).  Analogously, a non-negative
vector $z_{\rm o}$ is said to be an {\bf out-cover} of $w_{\rm o}$ if
$\delta _{z_{\rm o}}(v)\geq w_{\rm o}(v)$ holds for every node $v\in
V$.  
We say that a pair $(z_{\rm o},z_{\rm i})$ of non-negative vectors defined on $\olra A$
is a {\bf cover} of $(w_{\rm o},w_{\rm i})$ if $z_{\rm o}$ is an out-cover of $w_{\rm o}$ and $z_{\rm i}$ is an in-cover of $w_{\rm i}$. 
A cover $(z_{\rm o},z_{\rm i})$ {\add  of $(w_{\rm o},w_{\rm i})$} is a circular cover {\add of $(w_{\rm o},w_{\rm i})$} if $z = z_{\rm o} + z_{\rm i}$
is a circulation {\add in $\olra {D}$.}  
Let $c:= \chi_A$ denote the characteristic function of $A$ defined on
the arc-set $\olra A$, that is,

\eq c(e):= \begin{cases} 1 & \ \ \hbox{if}\ \ \ e\in A \cr 0 & \ \
\hbox{if}\ \ \ e\in \ola A. \label{(cdef)} \end{cases} \eeq

\noindent In what follows, the cost-function $c$ always refers to the
function defined in \eref{(cdef)}.  The $c$-cost of a circular cover $(z_{\rm o},z_{\rm i})$ {\add of $(w_{\rm o},w_{\rm i})$} 
is $cz_{\rm
o}+cz_{\rm i}$.

The main result of the paper (beside the algorithmic aspects) is the following min-max formula.
Recall that for a vector
$w$ on set $V$, $\widetilde w(X):=\sum \ [w(s):s\in X]$.

\begin{theorem} \label{minmax2} Let $w_{\rm o}:V\rightarrow {\bf R}_+$ and
$w_{\rm i}:V\rightarrow {\bf R}_+$ be weight-functions on the node-set
of digraph $D=(V,A)$.  Then

\eq \begin{cases} \hbox{ $\max \{\widetilde w_{\rm o}(Y_{\rm o})+
\widetilde w_{\rm i}(Y_{\rm i}):  (Y_{\rm o},Y_{\rm i})$ \ a
source-sink pair$\}$}\ \ & \cr \hbox{ $=$ }\ & \cr \hbox{ $\min
\{c(z_{\rm o}+z_{\rm i}):  \ (z_{\rm o},z_{\rm i})$ \ a circular cover
of $(w_{\rm o},w_{\rm i})$$\}$.}\ \cr \end{cases} \eeq If $(w_{\rm
o},w_{\rm i})$ is integer-valued, then the circular cover of
minimum $c$-cost
may be chosen integer-valued.  \end{theorem}

\subsection{A primal circulation problem and its dual tension problem}\label{subsec:circulation}

We associate a digraph $D\sp *=(V\sp *,A\sp *)$ with $D$ and investigate a
minimum cost feasible circulation problem on $D\sp *$ along with its
linear programming dual which is a maximum weight feasible tension
problem.  On one hand, an optimal integer-valued tension in $D\sp *$
shall define a maximum $(w_{\rm o},w_{\rm i})$-weight so-si pair.  On
the other hand, an optimal feasible circulation in $D\sp *$ shall
define a minimum $c$-cost circular cover $(z_{\rm o},z_{\rm i})$ of
$(w_{\rm o},w_u)$.

Let $V_{\rm i}$ and $V_{\rm o}$ be disjoint copies of $V$ and let
$V\sp *:=V_{\rm o}\cup V\cup V_{\rm i}$.  (For an intuition, we think
on these three sets to be drawn in three parallel horizontal lines so
that $V_{\rm o}$ is the lower, $V$ is the middle, and $V_{\rm i}$ is
the upper set, and their elements are positioned in such a way that
$v_{\rm o}, \ v, \ v_{\rm i}$ are in a vertical position for every
$v\in V$.)

Define a digraph $D\sp *=(V\sp *,A\sp *)$ as follows.  For every node
$v\in V$, let $v_{\rm i}v$ and $vv_{\rm o}$ be \lq vertical\rq \ arcs
of $D\sp *$ (which are directed downward).  Furthermore, for every arc
$uv\in \olra {A}$, let $uv$, $uv_{\rm i}$, and $u_{\rm o}v$ be arcs in
$D\sp *$.  (Therefore $D\sp *$ has $3\vert V\vert $ nodes and $2\vert
V\vert + 6\vert A\vert $ arcs, and the restriction of $D\sp *$ to $V$
is just $\olra{D}$. See Figure \ref{fig:D*}.)

\begin{figure}[ht]
\begin{center}
\begin{tikzpicture}[scale=1.2, transform shape]

  \pgfmathsetmacro{\d}{4}	
  \pgfmathsetmacro{\b}{2}
  \pgfmathsetmacro{\g}{1}
		
    \node[vertex, label=left:$u$] (u) at (0,0) {};
	\node[vertex, label=right:$v$] (v) at ($(u) + (\d, 0)$) {};
	\draw [-stealth, thick] (u) -- (v);

    \node[label=left:$V$] (V) at ($(u) + (\d * 2, 0)$) {};

    \node[] (b) at ($(u) + (-\g/2 , -\b/2)$) {};
    \node[] (j) at ($(u) + (\d * 2 +\g/2 , -\b/2)$) {}                  edge[dashed] (b);

	\node[vertex, label=left:$u_{\text{i}}$] (ui) at ($(u) + (0,-\b)$) {};
	\node[vertex, label=right:$v_{\text{i}}$] (vi) at ($(u) + (\d,-\b)$) {};

    \node[label=left:$V_{\text{i}}$] (Vi) at ($(u) + (\d * 2, -\b)$) {};

 	\node[vertex, label=left:$u$] (u') at ($(u) + (0,-\b * 2)$) {};
	\node[vertex, label=right:$v$] (v') at ($(u) + (\d,-\b * 2)$) {};

    \node[label=left:$V$] (V') at ($(u) + (\d * 2, -\b * 2)$) {};

 	\node[vertex, label=left:$u_{\text{o}}$] (uo) at ($(u) + (0,-\b * 3)$) {};
	\node[vertex, label=right:$v_{\text{o}}$] (vo) at ($(u) + (\d, - \b * 3)$) {};

    \node[label=left:$V_{\text{o}}$] (Vo) at ($(u) + (\d * 2, -\b * 3)$) {};

	\draw [-{stealth}, thick] (ui) -- (u');
	\draw [-{stealth}, thick] (u') -- (uo);
	\draw [-{stealth}, thick] (vi) -- (v');
	\draw [-{stealth}, thick] (v') -- (vo);

    \draw [-{stealth}, thick] (uo) edge [bend right, thick] (v');
	\draw [-{stealth}, thick] (u') edge [bend left, thick] (vi);
	\draw [-{stealth}, thick] (u') edge [bend right, thick] (v');

    \draw [-{stealth}, dotted, thick] (v') edge [bend right, dotted, thick] (ui);
	\draw [-{stealth}, dotted, thick] (vo) edge [bend left, dotted, thick] (u');
\draw[-{stealth}, dotted]{}   (v')     edge [bend right, dotted, thick] (u');

\end{tikzpicture}
\end{center}

\caption{Auxiliary digraph $D^*$ associated with $D$}
\label{fig:D*}
\end{figure}

On the arc-set of $D\sp *$, we define a lower bound function $f\sp
*:A\sp *\rightarrow {\bf R}\sp +$ and a cost-function $c\sp *:A\sp
*\rightarrow \{0,1\}$, as follows.

\eq f\sp *(e):= \begin{cases} w_{\rm i}(v) & \ \ \hbox{if}\ \ \
e=v_{\rm i}v \ \ (v\in V) \cr w_{\rm o}(v) & \ \ \hbox{if}\ \ \
e=vv_{\rm o} \ \ (v\in V) \cr 0 & \ \ \hbox{otherwise,}\
\label{(f*def)} \end{cases} \eeq

\eq c\sp *(e):= \begin{cases} 1 & \ \ \hbox{if}\ \ \ e=uv \ \ (uv\in
A) \cr 1 & \ \ \hbox{if}\ \ \ e=uv_{\rm i} \ \ (uv\in A) \cr 1 & \ \
\hbox{if}\ \ \ e=u_{\rm o}v \ \ (uv\in A) \cr 0 & \ \
\hbox{otherwise.}\ \label{(c*def)} \end{cases} \eeq

Let $Q\sp *$ denote the $(0,\pm 1)$-valued incidence matrix of $D\sp
*$, in which the columns correspond to the nodes and the rows
correspond to the arcs of $D\sp *$.  An entry of the row corresponding
to $e=uv$ is $+1$ in the column of $v$ and $-1$ in the column of $u$,
and 0 otherwise.  Consider the following (primal) polyhedron of
feasible circulations.

\eq R_{\rm pr}:= \{ x\sp * :  \ Q\sp *x\sp * = 0, \ x\sp * \geq f\sp *
\}.  \label{(Rprimal)} \eeq

Below we refer to the members of $R_{\rm pr}$ as {\bf feasible}
circulations.  
Consider the primal linear program

\eq {\rm OPT}_{\rm pr} :  = \min \{ c\sp *x\sp * :  Q\sp *x\sp * = 0,
\ x\sp * \geq f\sp * \} \label{(optprimal)}, \eeq

\noindent which is a cheapest feasible circulation problem.

The dual polyhedron belonging to \eref{(optprimal)} and the dual
linear program are the following

\eq R_{\rm du}:= \{(\pi \sp *,y\sp *):  \pi \sp *Q\sp * + y\sp * =
c\sp *, \ y\sp * \geq 0 \}, \label{(Rdual)} \eeq

\eq {\rm OPT}_{\rm du} := \max \{f\sp *y\sp *:  \pi \sp *Q\sp * + y\sp
* = c\sp *, \ y\sp * \geq 0\}, \label{(optdual)} \eeq

\noindent where $\pi \sp *=(\pi _{\rm o},\pi ,\pi _{\rm i})$ is a
function defined on $V\sp *$, while $y\sp *$ is defined on $A\sp *$.

The duality theorem of linear programming states that OPT$_{\rm pr} =$
OPT$_{\rm du}$.  As the describing incidence matrix $Q\sp *$ is known
to be totally unimodular (TU) and $c\sp *$ is $(0,1)$-valued, the dual
optimum is attained at an integral vector.  If $f\sp *$ is
integer-valued (that is, if $(w_{\rm o},w_{\rm i})$ is
integer-valued), then the primal optimum is also attained at an
integral vector.

The constraint equality $\pi \sp *Q\sp * + y\sp * = c\sp *$ in the
dual program requires for every arc $a\in A\sp *$ that

$$ \Delta _{\pi \sp *}(a) + y\sp *(a) = c\sp *(a).$$

In
particular, it follows that

\eq y\sp *(a) = - \Delta _{\pi \sp *}(a) \ \hbox{holds for each
vertical arc $a\in A\sp *$ of $D\sp *$. }\ \label{(ycsillag2)} \eeq

\noindent Since $y\sp *$ is non-negative, $\pi \sp *$ is a  $c\sp *$-feasible
potential, that is, $ \Delta _{\pi \sp *} (a)
\leq c\sp *(a)$ \ holds for each arc \ $a\in A\sp *$.  Therefore, if
we associate a vector $y\sp *:= c\sp * - \Delta _{\pi \sp *}$ (on $A\sp
*$) with a $c\sp *$-feasible potential $\pi \sp *$, the vector $(\pi \sp
*,y\sp *)$ obtained in this way is in the dual polyhedron
\eref{(Rdual)}.  (This means that the polyhedron of $c\sp
*$-feasible potentials is the projection of $R_{\rm du}$.)  
By applying (\ref{(ycsillag2)}), we get that the weight $M(\pi \sp *)$ of $\pi \sp *$ in the
dual linear program is

\eq M(\pi \sp *):= f\sp *y\sp * = \sum \ [-\Delta _{\pi \sp *}(a)f\sp *(a):  \
a \in A\sp * \ \hbox{a vertical arc}].  \label{(Mdef)} \eeq

\subsection{From feasible potentials in {\boldmath$D^*$} to so-si pairs in {\boldmath$D$} }

Consider an integer-valued $c\sp *$-feasible potential $\pi \sp *=(\pi
_{\rm o},\pi ,\pi _{\rm i})$ along with the vector $y\sp *$ for which
$(\pi \sp *,y\sp *)\in R_{\rm du}$.  We show how to associate a so-si
pair $(Y_{\rm o},Y_{\rm i})$ with $\pi \sp *$ for which

\eq \widetilde w_{\rm o}(Y_{\rm o})+ \widetilde w_{\rm i}(Y_{\rm i})=
M(\pi \sp *) \ \ (= f\sp *y\sp *)  \label{(yoyi)} \eeq

\noindent holds.

\begin{claim}\label{steady2} The potential $\pi $ (on $V$) occurring in $\pi
\sp *=(\pi _{\rm o},\pi ,\pi _{\rm i})$ is small-dropping on the
arc-set of $D$, that is, $\Delta _{\pi }(e)\in \{0,1\}$ holds for each
arc $e\in A$.  
\end{claim} 
\Proof What we prove is that $0\leq \Delta _{\pi }(e)\leq 1$.  On one
hand, clearly $$ \Delta _\pi (e) = \Delta _{\pi \sp *}(e) \leq c\sp
*(e)=1.$$ On the other hand, $\Delta _\pi (e) + \Delta _\pi (\ola{e})
=0$ holds for the arc $\ola e \in \ola A$ arising by reversing $e$,
and hence

$$ - \Delta _\pi (e) = \Delta _\pi (\ola e) = \Delta _{\pi \sp *}(\ola
e) \leq c\sp *(\ola e)=0,$$ that is, $\Delta _\pi (e)\geq 0$.  \FBOX

\begin{claim}\label{fuggoleges} For every node $v\in V$, the sum $y\sp
*(v_{\rm i}v)+y\sp *(vv_{\rm o})$ is either $0$ or $1$.  In
particular, the non-negativity of $y\sp *$ implies that $y\sp *(v_{\rm
i}v)\in \{0,1\}$ and $y\sp *(vv_{\rm o})\in \{0,1\}.$
\end{claim}

\Proof Assume first that there is an arc $vu\in A$ leaving $v$.  By
\eref{(ycsillag2)} we have

$$y\sp *(v_{\rm i}v) +y\sp *(vv_{\rm o}) = [\pi \sp *(v_{\rm i}) - \pi
\sp *(v)] + [\pi \sp *(v) - \pi \sp *(v_{\rm o})] = \pi \sp *(v_{\rm
i}) - \pi \sp *(v_{\rm o}) =$$ $$[\pi \sp *(v_{\rm i})-\pi \sp *(u)] +
[\pi \sp *(u) -\pi \sp *(v_{\rm o})] \leq c\sp *(uv_{\rm i}) + c\sp
*(v_{\rm o}u) = 0 +1 =1.$$

If no arc of $D$ leaves $v$, then there is an arc $uv\in A$ entering
$v$, for which the proof runs analogously:

$$y\sp *(v_{\rm i}v) +y\sp *(vv_{\rm o}) = [\pi \sp *(v_{\rm i}) - \pi
\sp *(v)] + [\pi \sp *(v) - \pi \sp *(v_{\rm o})] = \pi \sp *(v_{\rm
i}) - \pi \sp *(v_{\rm o}) =$$ $$[\pi \sp *(v_{\rm i})-\pi \sp *(u)] +
[\pi \sp *(u) -\pi \sp *(v_{\rm o})] \leq c\sp *(uv_{\rm i}) + c\sp
*(v_{\rm o}u) = 1 +0 =1.$$ \FBOX

\medskip

Define subsets $Y_{\rm o}$ and $Y_{\text{i}}$ of $V$, as follows.

\eq \begin{cases} Y_{\rm o}:  = \{v\in V:  y\sp *(vv_{\rm o}) = 1 \},
\cr Y_{\rm i}:  = \{v\in V:  y\sp *(v_{\rm i}v)=1\}.
\label{(YoYidef)} \end{cases} \eeq

\begin{lemma}\label{pisosi} The pair $(Y_{\rm o},Y_{\rm i})$ of subsets
defined in (\ref{(YoYidef)}) 
forms a source-sink pair for which \emeref{(yoyi)} holds.
\end{lemma}

\Proof Observe first that $Y_{\rm o}\cap Y_{\rm i}=\emptyset $ since
if there were a node $v\in Y_{\rm o}\cap Y_{\rm i}$, then $y\sp
*(v_{\rm i}v) = 1 = y\sp *(vv_{\rm o})$, contradicting Claim
\ref{fuggoleges}.  Furthermore, we have

$$\widetilde w_{\rm o}(Y_{\rm o}) = \sum \ [w_{\rm o}(v)y\sp *(vv_{\rm
o}) :  v\in V] = \sum \ [f\sp *(vv_{\rm o})y\sp *(vv_{\rm o}) :  v\in
V], $$ \noindent and

$$\widetilde w_{\rm i}(Y_{\rm i}) = \sum \ [w_{\rm i}(v)y\sp *(v_{\rm
i}v) :  v\in V] = \sum \ [f\sp *(v_{\rm i}v)y\sp *(v_{\rm i}v) :  v\in
V], $$ from which \eref{(yoyi)} follows immediately.

By Claim \ref{steady2}, $\pi$ is small-dropping on the arc-set of $D$.
Next, we show that $Y_{\rm o}$ is a source set in the dicut-equivalent reorientation $D'$ of $D$ defined by $\pi$ so that we reverse those arcs $uv$ of $A$ for which $\Delta_\pi(uv) = 1$.   
Consider an arbitrary node $v\in Y_{\rm o}$.  By Theorem 
\ref{so-si.char}, it is enough to show for every arc $vt\in A$ leaving $v$
is that $\pi \sp *(t) - \pi \sp *(v)=0$, and for every arc $sv\in A$
entering $v$ that $\pi \sp *(v) - \pi \sp *(s)=1$.

By \eref{(ycsillag2)} it follows for $v_{\rm o}\in V_{\rm o}$ that
$\pi \sp *(v_{\rm o})-\pi \sp *(v) = - y\sp *(v_{\rm o}v) = -1$, and
hence we have for each arc $vt\in A$ the following.

$$\pi \sp *(t)-\pi \sp *(v) = [\pi \sp *(t)-\pi \sp *(v_{\rm o})] +
[\pi \sp *(v_{\rm o})-\pi \sp *(v)] = $$ $$ [\pi \sp *(t)-\pi \sp
*(v_{\rm o})] -1 \leq c\sp *(v_{\rm o}t) -1 = 0,$$

\noindent that is, $\pi \sp *(t)-\pi \sp *(v) \leq 0$.  On the other
hand, as $tv\in \ola A$, we have $\pi \sp *(v)-\pi \sp *(t) \leq c\sp
*(tv) = 0$, from which $0\leq \pi \sp *(t) - \pi \sp *(v) \leq 0$,
that is, $\pi \sp *(t) - \pi \sp *(v)=0$ follows.

Similarly, for an arc $sv\in A$ we get:

$$\pi \sp *(s)-\pi \sp *(v) = [\pi \sp *(s)-\pi \sp *(v_{\rm o})] +
[\pi \sp *(v_{\rm o})-\pi \sp *(v)] = $$ $$ [\pi \sp *(s)-\pi \sp
*(v_{\rm o})] -1 \leq c\sp *(v_{\rm o}s) -1 = -1,$$

\noindent from which $\pi \sp *(v)-\pi \sp *(s) \geq 1$.  On the other
hand $\pi \sp *(v)-\pi \sp *(s)\leq c\sp *(cs)=1$, and hence $1\leq
\pi \sp *(v) - \pi \sp *(s) \leq 1$, that is, $\pi \sp *(s) - \pi \sp
*(v)=1$ follows.

An analogous argument shows that $Y_{\rm i}$ is a sink set in $D'$, completing the proof of the lemma.  \FBOX

\medskip It follows from Lemma \ref{pisosi} that:

\begin{eqnarray} \max \{\widetilde w_{\rm o}(Y_{\rm o})+\widetilde
w_{\rm i}(Y_{\rm i}) :  (Y_0,Y_{\rm i}) \ \hbox{is a so-si pair}\} & \geq
& \cr \max \{M(\pi \sp *):  \ \pi \sp * \ \ \hbox{is a $c\sp *$-feasible
potential}\} \ & = & {\rm OPT}_{\rm du}.  \label{(rezgeq)}
\end{eqnarray}

\subsection{From feasible circulations in {\boldmath$D\sp *$} to circular covers
of {\boldmath$(w_{\rm o}, w_{\rm i})$} } \label{fotetel.biz}

Consider an optimal solution $z^*$ of  \eref{(optprimal)}, which may be assumed to be integer-valued when $(w_{\rm o},w_{\rm i})$  is integer-valued.
Next, we associate a circular cover $(z_{\rm o},z_{\rm i})$ of
$(w_{\rm o},w_{\rm i})$ with  $z\sp *$ in such a way that $c(z_{\rm o}+z_{\rm i})
= c\sp *z\sp *$ (where $c=\chi_A$ is defined in \eref{(cdef)}).  Here
both $z_{\rm o}$ and $z_{\rm i}$ are non-negative vectors on the
arc-set $\olra{A}$, and their sum $z_{\rm o}+z_{\rm i}$ is a
circulation in $\olra D$.

We may assume that $z\sp *(uv)=0$ for every arc $uv\in A$ since if
$\alpha =z\sp *(uv)$ were positive
, then by
decreasing $z\sp *(uv)$ to $0$ and increasing both $z\sp *(uv_{\rm
i})$ and $z\sp *(v_{\rm i}v)$ by $\alpha $, we obtain a new
circulation which is also $f\sp *$-feasible and its $c\sp *$-cost is
unchanged.

Define functions $z_{\rm i}:\olra {A}\rightarrow {\bf R_+} $ and
$z_{\rm o}:\olra {A}\rightarrow {\bf R_+} $ as follows.

\eq z_{\rm o}(uv):= z\sp *(u_{\rm o}v) \ \hbox{and}\ \ z_{\rm i}(uv):=
z\sp *(uv_{\rm i}).  \label{(zozidef)} \eeq

Note that this $(z_{\rm o},z_{\rm i})$  is integer-valued  when $(w_{\rm o},w_{\rm i})$  is integer-valued. 
Observe that the pair $(z_{\rm o},z_{\rm i})$  may be viewed as a vector arising from $z^*$ by deleting the
arcs in $A$ and contracting the arcs $v_{\rm o}v$ and $v_{\rm i}v$ for every $v\in V$. 
Hence the sum $z:=z_{\rm o}+z_{\rm i}$ is a circulation in $\olra {D}$
for which $cz= c\sp *z\sp *$.  Moreover, we have for every node $u\in
V$ the following.

\eq \begin{cases} \delta _{z_{\rm o}}(u)=\delta _{z\sp *}(u_{\rm o}) =
z\sp *(uu_{\rm o}) \geq f\sp *(uu_{\rm o}) = w_{\rm o}(u) \cr \varrho
_{z_{\rm i}}(u)=\varrho _{z\sp *}(u_{\rm i}) = z\sp *(u_{\rm i}u) \geq
f\sp *(u_{\rm i}u) = w_{\rm i}(u), \end{cases} \label{(kibefed)} \eeq

\noindent from which it follows that $(z_{\rm o},z_{\rm i})$ is indeed
a circular cover of $(w_{\rm o},w_{\rm i})$.

This implies that:

\begin{eqnarray} \min \{c(x_{\rm o}+x_{\rm i}):  (x_{\rm o},x_{\rm i}) \ \hbox{a circular cover of}\ (w_{\rm o},w_{\rm i})\} \leq \cr c(z_{\rm o}+z_{\rm i}) \ =  \ c\sp *z\sp * \ =\ {\rm OPT}_{\rm pr}.\label{(kibeleq)}
\end{eqnarray}

\section{The proof of Theorem \ref{minmax2}} \label{mainproof}

\subsection{Proving {\boldmath $\max\leq \min$} }

\begin{claim} \label{cekvi} Let $D'=(V,A')$ be a dicut-equivalent
reorientation of $D$ and let $c':= \chi_{A'}$ denote the characteristic
function of $A'$ on $\olra A$.  Then $cz=c'z$ holds for every
circulation $z\geq 0$ (defined on digraph $\olra D$).  \end{claim}

\Proof Recall the notation $\widetilde c$.  Since a non-negative
circulation can be expressed as the non-negative linear combination of
one-way circuits it suffices to show that $\widetilde c(K) =
\widetilde c'(K)$ holds for every one-way circuit $K$ of $\olra D$.
But this will follow once we prove it for the case when $D'$ arises
from $D$ by reversing a single dicut.

Let this dicut be the set of arcs of $D$ entering a subset $Z\subset
V$ of nodes with no leaving arcs of $D$.  In this case $c'$ arises
from $c$ in such a way that the $c$-cost of $A$-arcs entering $Z$ is
changed from $1$ to $0$ while the $c$-cost of $\ola{A}$-arcs leaving $Z$
is changed from $0$ to $1.$ But this implies, as $\varrho _K(Z)=\delta
_K(Z)$, that the cost of one-way circuit $K$ does not change, that is,
$\widetilde c(K) = \widetilde c'(K)$.  \FBOX

\medskip

\Remark \label{rovbiz.cekvi}  \emph{A more
concise proof of Claim \ref{cekvi}, where the non-negativity of $z$ is
not used, is as follows.  By Claim \ref{vekvi1}, there is a potential
$\pi \in V\rightarrow {\bf Z} $ (which is small-dropping on $A$) for
which $A'$ arises from $A$ by reorienting those arcs $e\in A$ for
which $\Delta _\pi (e)=1$.  Now $c'=c-\Delta _\pi $ and $z\Delta _\pi
=0$ since the scalar product of a circulation and a tension is always
zero.  Therefore $c'z=cz$ follows indeed.}  $\bullet $ \eRe \medskip

The requested inequality $\max \leq \min$ immediately follows from the
following claim.

\begin{claim} \label{cDz} Let $(Y_{\rm o},Y_{\rm i})$ be a so-si pair in $D$,
 and let $(z_{\rm
o},z_{\rm i})$ be a circular cover of $(w_{\rm o},w_{\rm i})$.  Then
$cz\geq \widetilde w_{\rm o}(Y_{\rm o}) + \widetilde w_{\rm i}(Y_{\rm
i})$ holds for the circulation $z=z_{\rm o}+z_{\rm i}$.  \end{claim}

\Proof Consider the dicut-equivalent reorientation $D'$ of $D$ for
which the elements of $Y_{\rm o}$ are source nodes and the elements of
$Y_{\rm i}$ are sink nodes.  For a vector $x:\olra A\rightarrow {\bf
Z}$ and for a subset $Z\subset V$, let $\delta '_x(Z):=\sum \ [x(e):  e\in A', \ e$ leaves $Z]$ and
$\varrho '_x(Z):=\sum \ [x(e):  e\in A', \ e$ enters $Z$$]$.  By Claim
\ref{cekvi}, $cz=c'z \ (=c'z_{\rm o}+c'z_{\rm i})$.  But we have
$c'z_{\rm o} \geq \delta '_{z_{\rm o}}(Y_{\rm o})\geq \widetilde
w_{\rm o}(Y_{\rm o})$ and $c'z_{\rm i} \geq \varrho '_{z_{\rm
i}}(Y_{\rm i})\geq \widetilde w_{\rm i}(Y_{\rm i})$ from which $cz\geq
\widetilde w_{\rm o}(Y_{\rm o}) + \widetilde w_{\rm i}(Y_{\rm i})$
follows.  \FBOX \medskip

\subsection{Proving {\boldmath $\max \geq \min$} } \label{maxatleastmin}

By combining \eref{(kibeleq)} with inequality \eref{(rezgeq)}, we get
the following.

\eq \begin{cases} \hbox{$\max \{\widetilde w_{\rm o}(Y_{\rm o}) +
\widetilde w_{\rm i}(Y_{\rm i}):  (Y_{\rm o},Y_{\rm i})$ \ a so-si pair$\} \ \geq $}\ & \cr {\rm OPT}_{\rm du} \ = \ {\rm OPT}_{\rm pr} \
\geq & \cr \hbox{$ \min \{c(z_{\rm o}+z_{\rm i}):  (z_{\rm o},z_{\rm
i})$ a circular cover of $(w_{\rm o},w_{\rm i})\}$, }\
\end{cases} \label{(maxatleast)} \eeq

\medskip

\noindent and this is exactly the non-trivial inequality $\max \geq
\min$ of Theorem \ref{minmax2}.

Finally, observe that if $(w_{\rm o},w_{\rm i})$ is integer-valued,
then $f\sp *$ in the primal linear program \eref{(optprimal)} is also
integer-valued, and hence the optimal circulation can be chosen
integer-valued.  In this case, the pair $(z_{\rm o},z_{\rm i})$
associated with $z\sp *$ in \eref{(zozidef)} is an integer-valued circular
cover of $(w_{\rm o},w_{\rm i})$, and hence the proof of Theorem
\ref{minmax2} is complete.  \BB

\Corollary \label{optrez} 
Let $z\sp *$ be an optimal solution to the
primal linear program \emeref{(optprimal)}, and let $(z_{\rm o},z_{\rm
i})$ be the circular cover of $(w_{\rm o},w_{\rm i})$ associated with
$z\sp *$ in \emeref{(zozidef)}, which is integer-valued when $(w_{\rm
o},w_{\rm i})$ is integer-valued.  Let $(\pi \sp *,y\sp *)$ be an
optimal integral solution to the dual linear program
\emeref{(optdual)}, and let $(Y_0,Y_{\rm i})$ be the so-si pair
associated with $\pi \sp *$ in \emeref{(YoYidef)}.  Then $c(z_{\rm
o}+z_{\rm i}) = \widetilde w_{\rm o}(Y_0)+ \widetilde w_{\rm i}(Y_{\rm
i})$.  Moreover, $(z_{\rm o},z_{\rm i})$ is a minimum $c$-cost
circular cover of $(w_{\rm o},w_{\rm i})$, and $(Y_0,Y_{\rm i})$ is a
maximum $(w_{\rm o},w_{\rm i})$-weight so-si pair. \FBOX \eCo

\subsection{Algorithmic aspects} \label{algasp}
In the proof of
Theorem \ref{minmax2}, we introduced a network circulation problem on
a digraph $D\sp *=(V\sp *,A\sp *)$ along with its dual tension
problem (see Formulas (\ref{(optprimal)}) and (\ref{(optdual)})) where values ${\rm OPT}_{\rm pr}={\rm OPT}_{\rm du}$ give the optimal value for the double-weighted source-sink pair problem (\ref{def:WSSP}). These circulation and tension problems are special in the sense that the (primal) cost-function $c\sp *$ on $A\sp *$ is $(0,1)$-valued.  It is a
standard reduction to reformulate such a circulation problem as a
minimum cost network flow problem in which the cost-function is
$(0,1)$-valued (see e.g. Section 3.4.2. in \cite{FrankCO} on the equivalence of circulations and flows).  The min-cost flow algorithm of Ford and Fulkerson
\cite{Ford-Fulkerson} in this case is strongly polynomial, provided
that a strongly polynomial subroutine is used for the intermediate
maximum flow computations (for example, the classic shortest
augmenting path type algorithms of Dinitz or Edmonds and Karp  or
the push-relabel algorithm of Goldberg and Tarjan).

Therefore, an optimal solution to the primal and to the dual linear
programs \eref{(optprimal)} and \eref{(optdual)} can be computed in
strongly polynomial time, where the dual optimal vector is
integer-valued and the primal optimal vector is integer-valued when
$(w_{\rm o},w_{\rm i})$ is integer-valued.  By Corollary \ref{optrez},
a maximum $(w_{\rm o},w_{\rm i})$-weight so-si pair occurring in
Theorem \ref{minmax2} as well as a minimum $c$-cost circular cover of
$(w_{\rm o},w_{\rm i})$ can be computed in strongly polynomial time.

\subsection{Consequences}

Theorem 5.1 in \cite{FrankJ66} provided a min-max formula for the
maximum $w$-weight of a sink-stable set of a digraph $D=(V,A)$ for an integer-valued weight-function $w$.  Here we show how this formula follows from Theorem \ref{minmax2}.

We say that a family {\add (with the agreement that such a
family is allowed to contain repeated copies of a given member)} of one-way circuits of digraph $\olra D$ \textbf{covers} $w$
if every node $u\in V$ belongs to at least $w(u)$ circuits.  The {\bf
$A$-value} of a one-way circuit $K$ of digraph $\olra D$ is $\vert
K\cap A\vert $. In particular, $\{e,\ola e\}$ is a two-element one-way
circuit of $\olra D$ for every arc $e\in A$, whose $A$-value is $1$, thus a cover always exists for every weighting $w$ if $D$ is weakly connected and has at least two nodes.  
With this terminology, Theorem 5.1 in \cite{FrankJ66} is equivalent to
the following.\footnote{A slightly different notion was used for a cover in \cite{FrankJ66}. Here we define a cover as a family of one-way circuits of $\olra D$, whereas a cover in \cite{FrankJ66} consisted of circles and arcs of $D$.}

\begin{theorem} \label{mimax.wClar} Let $w$ be a non-negative, integer-valued weight-function on the arc-set of a
digraph $D=(V,A)$. Then the maximum $w$-weight of a sink-stable
set of $D$ is equal to the minimum sum of $A$-values
of a family of one-way circuits of digraph $\olra D$ that cover $w$.
\end{theorem}

\Proof Let $w_{\rm o}:\equiv 0$ and $w_{\rm i}:=w$.  
For a so-si pair $(Y_{\rm o},Y_{\rm i})$, the set $Y_{\rm i}$ is a sink-stable set for which $\widetilde w_{\rm o}(Y_{\rm o}) + \widetilde w_{\rm i}(Y_{\rm i}) = \widetilde w(Y_{\rm i})$. 
For a sink-stable set $Y$, 
pair $(\emptyset,Y)$ is a so-si pair for which  $\widetilde w_{\rm o}(\emptyset) + \widetilde w_{\rm i}(Y) = \widetilde w(Y)$. 
Thus the maximum $(w_{\rm o},w_{\rm i})$-weight of a so-si pair is equal to the maximum $w$-weight of a sink set.
By Theorem \ref{minmax2}, this latter maximum equals the minimum of $c(z_{\rm o}+z_{\rm i})$ over (non-negative, integer-valued) circular covers $(z_{\rm o},z_{\rm i})$ of $(0,w)$. 
We show that this minimum equals the  minimum sum of $A$-values
of a family of one-way circuits of digraph $\olra D$ that cover $w$.

First, let $(z_{\rm o},z_{\rm i})$ be a minimum cost circular cover of $(0,w)$.
Then $z=z_{\rm o}+z_{\rm i}$ is a circulation in $\olra D$ such that $(0,z)$ is a circular cover of $(0,w)$ and $c(z_{\rm o}+z_{\rm i})=cz$.  But $z$ can be obtained as a non-negative integral combination of one-way circuits of $\olra D$, and these coefficients define a family of one-way circuits. The property that $(0,z)$ covers $(0,w)$ means that every node $u$ belongs to at least $w(u)$ one-way circuits {\add in the family}. That is, this family of circuits covers $w$ and the sum of $A$-values of the circuits is $cz$.

Second, if 
a family of one-way circuits of digraph $\olra D$ covers $w$,
then the sum of the incidence vectors of the circuits in the family gives a circulation $z$ such that $(0,z)$ is a circular cover of $(0,w)$, and cost $c(0+z)$ equals the sum of $A$-values of the circuits in the family.
\FBOX

\medskip
\Remark
\emph{Theorem 5.1 of [7] (and hence the present
Theorem \ref{mimax.wClar}) does not imply directly the min-max theorem of Abeledo
and Atkinson (Theorem 8 of \cite{Abeledo-Atkinson4}) on the maximum cardinality of a Clar-set of bounded faces
of $G=(S,T;E)$.  Instead, a sharpening (Theorem 5.4 of [7]) of a special case of Theorem 5.1 of [7] (when there exists a so-called "thin" subset of edges) was derived in [7] from which the Abeledo-Atkinson theorem followed immediately.  Therefore the present Theorem \ref{mimax.wClar} can imply only a weaker version of the Abeledo-Atkinson theorem stating that, in a perfectly matchable plane bipartite graph, the Clar number is equal to the minimum total value of (minimal) cuts intersecting all bounded faces, where the value of a
cut $[Z, \ol Z]$ is the in-degree of $Z\subset S\cup Z$ in the digraph
${\ora G_M}$ arising from $G$ by orienting the elements of an
arbitrary perfect matching $M$ toward $S$ while the other edges of $G$ toward $T$.  The content of the original theorem of Abeledo-Atkinson is that the minimum can actually be 
attained on a family of such cuts $[Z, \ol Z]$ for which no element of $\ora M$ leaves $Z$ and every arc in ${\ora G}_M$ entering $Z$ belongs to $\ora M$. Finally, we note (without going into details) that it is not difficult to obtain this sharpened version form the weaker one.} 
 $\bullet $ 
\eRe

Next, we show a min-max formula concerning maximum weight resonant
sets.  
Recall the definition of a circular cover $(z_{\rm o},z_{\rm i})$ {\add  of $(w_{\rm o},w_{\rm i})$} introduced in Section \ref{subsec:MT} and the cost-function $c$ on arc-set $\olra A$
given in \eref{(cdef)}.  In the special case when $w:=w_{\rm o}=w_{\rm i}$,
Theorem \ref{minmax2} is as follows.

\begin{theorem} \label{minmax2b} Let $w:V\rightarrow {\bf R}_+$ be a
non-negative weight-function on the node-set of digraph $D=(V,A)$.
Then

\eq \begin{cases} \hbox{ $\max \{\widetilde w(Y) :  Y\subseteq V$ \
a resonant set$\}$}\ & \cr \hbox{ $=$ }\ & \cr \hbox{ $\min \{c(z_{\rm
o}+z_{\rm i}):  \ (z_{\rm o},z_{\rm i})$ a circular cover of
$(w,w)\}$.}\ \cr \end{cases} \eeq When $w$ is
integer-valued, a circular cover of minimum $c$-cost can be chosen
integer-valued.  \FBOX \end{theorem}

We say for a non-negative integer vector $z_{\rm i}$ on arc-set
$\olra{A}$ that it is an {\bf in-cover} of a subset $U\subseteq V$ of
nodes if $\varrho _{z_{\rm i}}(v)\geq 1$ holds for each node $v\in U$.
A non-negative integer vector $z_{\rm o}$ defined on arc-set
$\olra{A}$ is an {\bf out-cover} of $U$ if
$\delta _{z_{\rm o}}(v)\geq 1$ holds for every node $v\in U$.
Let $U_{\rm o}\subseteq V$ and $U_{\rm i}\subseteq V$ be two disjoint
subsets.  We say that the pair $(z_{\rm o},z_{\rm i})$ of vectors is a {\bf cover}
of the pair $(U_{\rm o},U_{\rm i})$ of sets if $z_{\rm o}$ is an
out-cover of $U_{\rm o}$ and $z_{\rm i}$ is an in-cover of $U_{\rm
i}$.   

By applying Theorem \ref{minmax2} to
weight-functions $w_{\rm o}:=\chi_{U_{\rm o}}$ and $w_{\rm
i}:=\chi_{U_{\rm i}}$, we obtain the following.

\begin{theorem} \label{minmaxUoUi} In a digraph $D=(V,A)$, let $U_{\rm o}
\subseteq V$ and $U_{\rm i}\subseteq V$ be two disjoint subsets of
nodes.  Then 
\eq \begin{cases} \hbox{ $\max \{\vert Y_{\rm o}\vert +\vert
Y_{\rm i}\vert :  (Y_{\rm o},Y_{\rm i})$ \ a so-si pair, \ $Y_{\rm
o}\subseteq U_{\rm o}, \ Y_{\rm i}\subseteq U_{\rm i}\}$}\ \ & \cr \hbox{ $=$ }\ & \cr \hbox{ $\min
\{c(z_{\rm o}+z_{\rm i}):  \ (z_{\rm o},z_{\rm i})$ \ a circular cover
of $(U_{\rm o},U_{\rm i})$$\}$.}\ \cr \end{cases} \eeq
\FBOX \end{theorem}

\medskip

In the special case when $U:=U_{\rm
o}=U_{\rm i}$, this reduces to the the following.

\Corollary \label{cor:U}
In a digraph $D=(V,A)$, the
maximum cardinality of a resonant subset of a given subset $U\subseteq
V$ is equal to the minimum $c$-cost of a circulation $z$ (in $\olra
D$) $(*)$ which arises as the sum of non-negative integer-valued vectors
$z_{\rm o}$ and $z_{\rm i}$ for which $\varrho _{z_{\rm o}}(u) \geq 1 $ and $\delta
_{z_{\rm i}}(u) \geq 1$ hold for each node $u\in U$ (that is, $(z_{\rm o},z_{\rm i})$ is a circular pair  covering the pair $(U,U)$). In particular, $U$ is
a resonant set if and only if $cz\geq \vert U\vert $ holds for every
circulation $z$ (in $\olra D$) satisfying $(*)$.  \FBOX \eCo

\medskip

\begin{figure}[ht]
\begin{center}
\begin{tikzpicture}[scale=1.3, transform shape]

  \pgfmathsetmacro{\d}{4}	
  \pgfmathsetmacro{\b}{1}
  \pgfmathsetmacro{\g}{1}
		
    \node[] (u) at (0,0) {};
	
    \node[vertex, label=above:$a_1$, fill=black] (a1) at ($(u) + (\b, 0)$) {};
    \node[vertex, label=left:$a_2$] (a2) at ($(u) + (0, -\b)$) {};
    \node[vertex, label=below:$a_3$] (a3) at ($(u) + (\b, -\b * 2)$) {};

    \node[vertex, label=above:$x$, fill=black] (x) at ($(u) + (\b * 2, -\b)$) {};

    \node[vertex, label=above:$b_1$, fill=black] (b1) at ($(u) + (\b * 3,0)$) {};
    \node[vertex, label=right:$b_2$] (b2) at ($(u) + (\b * 4,-\b)$) {};
    \node[vertex, label=below:$b_3$] (b3) at ($(u) + (\b * 3,-\b * 2)$) {};

	\draw [-stealth, thick] (x) -- (a1);
	\draw [-stealth, thick] (a2) -- (a1);
	\draw [-stealth, thick] (a3) -- (a2);
	\draw [-stealth, thick] (x) -- (a3);
	\draw [-stealth, thick] (b1) -- (x);
	\draw [-stealth, thick] (b1) -- (b2);
	\draw [-stealth, thick] (b2) -- (b3);
	\draw [-stealth, thick] (b3) -- (x);

\end{tikzpicture}
\end{center}
\caption{A digraph $D=(V,A)$ where $U = \{a_1, b_1, x\}$ is not resonant}
\label{fig:example}
\end{figure}

As an example, consider the digraph in Figure \ref{fig:example}. It is not difficult to check that the four-element set $\{a_2,a_3,b_1,b_2\}$ is resonant, on the other hand  we claim that set
$U:=\{a_1,b_1,x\}$ is not.
To show this, define $z_{\rm o}$ and $z_{\rm i}$ on the arc-set $\olra A$ as follows.
Let $z_{\rm o}$ be $1$ on the arcs $a_1a_2, a_2a_3, b_1x, xb_3$ while $0$
otherwise, and let $z_{\rm i}$ be $1$ on the arcs $b_3b_2, b_2b_1, a_3x,
xa_1$ while $0$ otherwise.  Then $z_{\rm i}$ is an in-cover of $U$ and $z_{\rm o}$
is an out-cover of $U$ for which $z:=z_{\rm i}+z_{\rm o}$ is a $(0,1)$-valued
circulation.  Furthermore $z_{\rm o}$ is $1$ only on one element of $A$
(namely, on $b_1x$) and $z_{\rm i}$ is $1$ only on one element of $A$
(namely, on $xa_1$), that is, $cz_{\rm o}=1$ and $cz_{\rm i}=1$, and hence $cz=
c(z_{\rm o}+z_{\rm i})=2$, implying that the 3-element node-set $\{a_1,b_1,x\}$
cannot indeed be a resonant set of $D$.

\Remark \label{kornemeleg} \emph{ There is a somewhat surprising discrepancy between the characterization of Clar-sets and resonant
sets. 
Theorem 4.1 in \cite{FrankJ66} states roughly that a stable set $U\subseteq V$ is sink-stable if and only if a certain inequality holds for every circuit of digraph $D$. This result was extended by Theorem \ref{so-si.char} in the present work where the equivalence of Parts (A) and (C) indicated that a pair $(Y_{\rm o},Y_{\rm i})$ of node-sets is a source-sink pair if and only if a certain inequality holds for
every circuit of $D$. On the other hand, the characterization of resonant sets in the second part of Corollary
\ref{cor:U} requires an inequality not only for circuits but for certain integer-valued circulations, as
well. The question naturally arises whether this characterization can be simplified in such a
way that a certain inequality should be required only for circuits and not for more general
circulations.
The example in Figure \ref{fig:example} demonstrates that such a simplification is not
possible. Here $D$ has exactly two circuits, and $U:=\{a_1,b_1,x\}$ is a set which is clearly resonant
when restricted to any of these two circuits, but $U$ is not a resonant set of $D$ as we pointed out above. Therefore no
characterization of resonant sets  may exist that requires an inequality only for single circuits.}
 $\bullet $ 
\eRe 

\medskip

Finally we remark that Theorem \ref{minmax2} can be applied to manage
further special cases.  For example, the one when we are interested in
finding a maximum weight resonant set including a specified subset.
Or more generally, one may be interested in finding a (maximum weight)
so-si pair $(Y_{\rm i},Y_{\rm o})$ for which $U_{\rm o}'\subseteq
Y_{\rm o}\subseteq U_{\rm o}$ and $U_{\rm i}'\subseteq Y_{\rm
i}\subseteq U_{\rm i}$.  By defining $w_{\rm o}$ ($w_{\rm i}$) to be
an appropriately large number on the elements of $U_{\rm o}'$ \
($U_{\rm i}'$), and zero on the elements of $V- U_{\rm o}$ \
($V-U_{\rm i}$), Theorem \ref{minmax2} and the algorithmic approach
described above can be applied.

\subsection{Conclusion}

In \cite{FrankJ66}, the notion of Clar-sets of a perfectly matchable plane bipartite graph was extended
to the one of sink-stable sets of a general digraph, along with a network flow approach for
developing min-max formulas and purely combinatorial algorithms for various generalizations,
including the one of finding a maximum weight Clar-set.
Fries-set is another important notion in mathematical chemistry, but unlike the situation
with Clar-sets, the literature does not seem to include any strongly polynomial purely combinatorial algorithm for computing the Fries number, the largest cardinality of a Fries-set.
In the present paper, we did provide such an algorithm, actually in a much more general
setting which may be considered as a common framework for Clar-sets and for Fries-sets.
Namely, we introduced the notion of a source-sink pair of node-sets as well as a resonant
node-set in a digraph, and developed a network flow approach that gave rise to min-max formulas and to purely combinatorial strongly polynomial algorithms for computing both a maximum double-weight source-sink pair and a maximum weight resonant subset. The algorithm also provides
an optimal solution to the dual problem formulated in the min-max theorem.

\subsection*{Acknowledgement}

We thank the anonymous reviewer for his/her careful reading of our manuscript and for the many
insightful comments and suggestions.
This research was supported by the ÚNKP-22-5 New National Excellence Program of the Ministry for Culture and Innovation from the source of the National Research, Development and Innovation Fund. This research has been implemented with the support provided by the Ministry of Innovation and Technology of Hungary from the National Research, Development and Innovation Fund, financed under the  ELTE TKP 2021-NKTA-62 funding scheme, and supported by the János Bolyai Research Scholarship of the Hungarian Academy of Science.

\end{document}